\documentclass[12pt]{amsart}

\pdfoutput=1

\usepackage[hyperindex,breaklinks]{hyperref}
\PassOptionsToPackage{unicode}{hyperref}
\PassOptionsToPackage{naturalnames}{hyperref}
\usepackage[a4paper]{geometry}
\usepackage{amsmath,amssymb,amsthm,amsfonts,color}
\usepackage{comment}

\usepackage{stmaryrd}
\usepackage{tikz}
\usepackage{bm}

\usepackage[shortlabels]{enumitem}
\usepackage{refcount}

\theoremstyle{definition}
\newtheorem{Def}{Definition}
\newtheorem{Fact}[Def]{Fact}
\newtheorem{Thm}[Def]{Theorem}
\newtheorem{Lem}[Def]{Lemma}

\newtheorem{Rem}[Def]{Remark}

\newtheorem{Ex}[Def]{Example}

\newtheorem{Claim}{Claim}
\newtheorem{Q}{Question}

\newtheorem{Thm?}[Def]{Theorem?}

\DeclareMathOperator{\GL}{\text{L}}
\DeclareMathOperator{\St}{\text{St}}

\newcommand{\Ord}{\ensuremath{\mathrm{Ord}}}
\newcommand{\ZFC}{\ensuremath{\mathsf{ZFC}}}

\newcommand{\FA}{\ensuremath{\text{{\sf FA}}}}
\newcommand{\MM}{\ensuremath{\text{{\sf MM}}}} 

\newcommand{\Qp}[1]{\left\llbracket #1 \right\rrbracket}
\newcommand{\bool}[1]{\mathsf{#1}}
\newcommand{\B}{\bool{B}}
\newcommand{\C}{\bool{C}}
\newcommand{\SSP}{\ensuremath{\text{{\sf SSP}}}}
\newcommand{\U}{\mathsf{U}}

\newcommand{\rstrct}{\upharpoonright}

\begin{document}


\title[Universally Baire sets in $2^{\kappa}$]{Universally Baire sets in $2^{\kappa}$}
%\title{Universally Baire sets in $2^{\kappa}$}
%\author{Daisuke Ikegami and Matteo Viale}
\date{\today}

\author[D.\ Ikegami]{Daisuke Ikegami}
\address[D.\ Ikegami]{Institute of Logic and Cognition, Department of Philosophy, Sun Yat-sen University, Xichang Hall 602, 135 Xingang west street, Guangzhou, 510275 CHINA}

\email[D.\ Ikegami]{\href{mailto:ikegami@mail.sysu.edu.cn}{ikegami@mail.sysu.edu.cn}}

\author[M.\ Viale]{Matteo Viale}
\address[M.\ Viale]{Dipartimento di Matematica \lq\lq Giuseppe Peano", Universit\`{a} di Torino, Via Carlo Alberto 10, 10123 Turin, Italy}

\email[M.\ Viale]{\href{mailto:matteo.viale@unito.it}{matteo.viale@unito.it}}

\thanks{The first author would like to thank the Japan Society for the Promotion of Science (JSPS) for its generous support through the grant with JSPS KAKENHI Grant Number 15K17586. The second author acknowledges support from the projects \emph{PRIN 2017-2017NWTM8R Mathematical Logic: models, sets, computability, PRIN 2022 Models, sets and classifications, prot. 2022TECZJA}, and from GNSAGA.}

\begin{abstract}
We generalize the basic theory of universally Baire sets of $2^\omega$ to a theory of universally Baire subsets of $2^\kappa$. We show that the fundamental characterizations of the property of being universally Baire have natural generalizations that can be formulated also for subsets of $2^\kappa$, in particular we provide four equivalent uniform definitions in the parameter $\kappa$ (for $\kappa$ an infinite cardinal) characterizing for each such $\kappa$ the class of universally Baire subsets of $2^\kappa$. For $\kappa=\omega$, these definitions bring us back to the original notion of universally Baire sets of reals given by Feng, Magidor and Woodin~\cite{MR1233821}. 
\end{abstract}

\maketitle

\section{Introduction}

Universally Baire sets of reals have played a prominent role in the analysis of determinacy axioms and to gauge the complexity of determinacy hypothesis providing a fine scale of complexity (by means of the Wadge rank) of different determinacy axioms. Furthermore the universally Baire sets 
were shown to define the largest class $\Gamma\subseteq \wp (2^\omega)$ whose theory is generically invariant in the presence of large cardinals:
for example one of Woodin's main results shows that letting $\mathsf{uB}$ denote the class of universally Baire subsets of $2^\omega$ in the presence of a supercompact cardinal $\kappa$ and class many Woodin cardinals which are limits of Woodin cardinals, after collapsing $2^{\kappa}$ to be countable by any set forcing, we have that the theory of the model $\GL(\Ord^\omega,\mathsf{uB})$ with real parameters cannot be changed by means of set sized forcing.

In this paper, we generalize the basic theory of universally Baire sets of $2^\omega$ to a theory of universally Baire subsets of $2^\kappa$ and start to analyze the properties of these sets. We give some interesting initial results but nonetheless many problems remain open, thus we expect that further researches on this topic can be pursued in several distinct directions.

We show that the fundamental characterizations of the property of being universally Baire have natural generalizations that can be formulated also for subsets of $2^\kappa$, in particular we provide four equivalent uniform definitions in the parameter $\kappa$ (for $\kappa$ an infinite cardinal) characterizing for each such $\kappa$ the class of $\kappa$-universally Baire subsets of $2^\kappa$. For $\kappa=\omega$, these definitions bring us back to the original notion of universally Baire sets of reals given by Feng, Magidor and Woodin~\cite{MR1233821}. 

Several lines of further investigation remain open, among many: can our results, which can be considered a form of descriptive set theory on the space $2^\kappa$, relate to the recent wave of interest on the possible generalizations to the space $2^\kappa$ of the classical descriptive set theory on $2^\omega$? 
Is $\wp (\omega_1)^{\#}$ an universally Baire subset of $2^{\omega_1}$ in models of $\mathsf{MM}^{+++}$ with respect to stationary set preserving complete Boolean algebras which preserve $\mathsf{MM}^{+++}$?
Can we define a sensible notion of Wadge rank for universally Baire subsets of $2^{\kappa}$? If so what are its properties?
In the final section, we list in some details these and other related open problems.

%To be filled in.

\section{Basic definitions and facts}\label{sec:basics}

From now on, we work in $\ZFC$ unless clearly specified. 

We assume that the reader is familiar with the basic theory of forcing (given in e.g.,~\cite{Jech} and \cite{VIABOOKFORCING}), basics of descriptive set theory (given in e.g., ~\cite{MR1321597}), and the basics of stationary tower forcings given in~\cite{MR2069032}.
%%% ADD PREREQUISITES AND BASIC REFERENCES

Throughout this paper, $\kappa$ will be an infinite cardinal, $\B$ will be a complete Boolean algebra, and $X$ will be a topological space.

Given an infinite cardinal $\kappa$, let $2^{\kappa } = \big\{ f \mid f \colon \kappa \to 2 = \{0,1\} \big\}$. From now on, we identify subsets of $\kappa$ with elements of $2^{\kappa}$ using their characteristic functions. We work with the {\it product} topology on $2^{\kappa}$, i.e., the basic open sets in $2^{\kappa}$ are of the form $[s] = \{ x \in 2^{\kappa} \mid x \supseteq s\}$ for some $s \colon \text{dom} (s) \to 2$ which is a finite partial function from $\kappa$ to $2$. %%%% FOOTNOTE!!!  
 It is easy to see that each $[s]$ is clopen in $2^{\kappa}$ and hence $2^{\kappa}$ is zero-dimensional. It is well-known that the space $2^{\kappa}$ is compact and Hausdorff. For a set $Y$, we write $\text{Fn}(\kappa , Y)$ for the collection of finite partial functions from $\kappa$ to $Y$. 

Given a complete Boolean algebra $\B$, we often work with the {\it Stone space} of $\B$ (denoted by $\text{St}(\B)$), that is the collection of ultrafilters on $\B$ topologized by the sets of the form $O_b  = \{ G \in \text{St}(\B) \mid b \in G\}$ for some $b \in \B$. 

Let us mention some basic properties of $\text{St}(\B)$. It is easy to see that each basic open set $O_b$ is clopen in $\text{St}(\B)$ and hence $\text{St}(\B)$ is zero-dimensional. It is well-known that $\text{St}(\B)$ is a compact Hausdorff space.  It is also a well-known fact that for any clopen subset $C$ of $\text{St}(\B)$, there is a unique $b \in \B$ such that $C = O_b$. Using this fact, for each clopen subset $C$ of $\text{St}(\B)$, let us write $b_C$ for the unique element of $\B$ with $O_{(b_C)} = C$. 

%We now discuss lifting up clopen subsets of $\text{St}(\B)$ and continuous functions from $\text{St}(\B)$ to $2^{\kappa}$ in a generic extension of $V$ via $\B$:
%\begin{Def}\label{lift-up}
%${}$
%
%\begin{enumerate}
%\item Let $C$ be a clopen subset of $\text{St}(\B)$ and $G$ be a $\B$-generic filter over $V$. Let $b_C$ be the unique element of $\B$ such that $C = O_{b_C}$. Then we define $C^{V[G]}$ to be the clopen set $O_{b_C}^{V[G]}$ in $V[G]$.
%
%\item Let $f$ be a continuous function from $\text{St}(G)$ to $2^{\kappa}$ and let $G$ be a $\B$-generic filter over $V$. We define the continuous function $f^{V[G]} \colon \text{St}(\B) \to 2^{\kappa}$ in $V[G]$ as follows: For all ultrafilters $H$ on $\B$ in $V[G]$ and all $\alpha < \kappa$,
%\begin{align*}
%f^{V[G]} (H) (\alpha) = 1 \iff H \in C_{\alpha}^{V[G]},
%\end{align*}
%where $C_{\alpha} = \{ G' \in \text{St}(\B) \mid f(G') (\alpha) =1 \}$. 
%\end{enumerate}
%\end{Def}
%
%It is easy to see that $C^{V[G]} \cap V = C$, that $f^{V[G]}$ is continuous in $V[G]$, and that $f^{V[G]} \upharpoonright V = f$.

A $\B$-name $\tau$ is called a {\it $\B$-name for a subset of $\kappa$} if $\Vdash_\B \lq\lq \tau \subseteq \check{\kappa}"$. We now discuss a fundamental correspondence between continuous functions from $\text{St}(\B)$ to $2^{\kappa}$, and $\B$-names for a subset of $\kappa$:
\begin{Def}\label{translations name cont}

${}$
\begin{enumerate}
\item For a continuous function $f$ from $\text{St}(\B)$ to $2^{\kappa}$, let $\tau_f = \{ (\check{\alpha} , b_{C_{\alpha}}) \mid \alpha < \kappa \}$, where $C_{\alpha}$ is the clopen subset of $\text{St}(\B)$ such that $C_{\alpha} = \{ G \in \text{St}(\B) \mid f(G) (\alpha) = 1\}$, and $b_{C_{\alpha}}$ is the unique element $b$ of $\B$ with $O_b = C_{\alpha}$. 

\item For a $\B$-name $\tau$ for a subset of $\kappa$, let $f_{\tau} \colon \text{St}(\B) \to 2$ be such that 
\begin{align*}
f_{\tau} (G) (\alpha) = 1 \iff  \Qp{\check{\alpha} \in \tau}_\B  \in G %%%%% BOOLEAN VALUE
\end{align*}
for all $\alpha < \kappa$ and $G \in \text{St}(\B)$, where $\Qp{\check{\alpha} \in \tau}_\B$ is the Boolean value of the statement \lq\lq $\check{\alpha} \in \tau$".
\end{enumerate}
\end{Def}

\begin{Lem}\label{cont-name}
Let $f$ be a continuous function from $\text{St}(\B)$ to $2^{\kappa}$ and let $\tau$ be a $\B$-name for a subset of $\kappa$.
\begin{enumerate}[(1)]
\item $\tau_f$ is a $\B$-name for a subset of $\kappa$.\label{enum-lem:cont-name1}

\item $f_{\tau}$ is a continuous function from $\text{St}(\B)$ to $2^{\kappa}$.\label{enum-lem:cont-name2}

\item $f_{(\tau_f)} = f$.\label{enum-lem:cont-name3}

\item $\Vdash_\B \lq\lq \tau_{(f_{\tau})} = \tau "$.\label{enum-lem:cont-name4}
\end{enumerate}
\end{Lem}

\begin{proof}
${}$

\begin{description}
\item[\ref{enum-lem:cont-name1}] Clear by the definition of $\tau_f$.

\item[\ref{enum-lem:cont-name2}] We show that $f_{\tau}$ is continuous. It is enough to show that for all $\alpha < \kappa$, the set $A_{\alpha} = \{ G \in \text{St}(\B) \mid f_{\tau} (G) (\alpha) = 1 \}$ is clopen in $\text{St}(\B)$. But by the definition of $f_{\tau}$, $f_{\tau} (G) (\alpha) = 1 \iff \Qp{\check{\alpha} \in \tau}_\B \in G$. Hence, setting $b_{\alpha} = \Qp{\check{\alpha} \in \tau}_\B$, $A_{\alpha} = O_{b_{\alpha}}$, which is clearly clopen in $\text{St}(\B)$. 

Similarly, one can verify that for all $G\in \St (\B)$ and $\alpha < \kappa$, we have $f_{\tau} (G) (\alpha) = 0 \iff \Qp{\check{\alpha} \notin \tau}_\B \in G$. So for all $G\in \St (\B)$, $f_{\tau} (G)$ is a function from $\kappa$ to $2$. Hence $f_{\tau}$ is a function from $\St (\B)$ to $2^{\kappa}$.  

\item[\ref{enum-lem:cont-name3}] We show that $f_{(\tau_f)} (G) (\alpha) = 1 \iff f(G) (\alpha)  =1$ for all $G \in \text{St}(\B)$ and $\alpha < \kappa$. In fact, we have the following equivalences:
\begin{align*}
f_{(\tau_f)} (G) (\alpha) = 1 \iff & \Qp{\check{\alpha} \in \tau_f}_\B \in G\\
				\iff   & b_{C_{\alpha}} \in G\\
				\iff   & G \in C_{\alpha}\\
				\iff   & f(G)(\alpha) = 1,
\end{align*}
where $C_{\alpha} = \{ H \in \text{St}(\B) \mid f(H) (\alpha) =1\}$. 

The first equivalence above follows from the definition of $f_{(\tau_f)}$. The third equivalence above follows from the property of $b_{C_{\alpha}}$, i.e., $O_{(b_{C_{\alpha}})} = C_{\alpha}$.\footnote{The notation, property, and uniqueness of $b_{C_{\alpha}}$ was discussed in the 5th paragraph of this section, i.e., Section~\ref{sec:basics}.} The fourth equivalence above follows from the definition of $C_{\alpha}$.

We will show the second equivalence above. In fact, we prove that $\Qp{\check{\alpha} \in \tau_f}_\B = b_{C_{\alpha}}$. Suppose that they are different. Then one can take a $\B$-generic filter $H$ over $V$ such that $\Qp{\check{\alpha} \in \tau_f}_\B$ is in $H$, but $b_{C_{\alpha}}$ is {\it not} in $H$. However, by the genericity of $H$, $\Qp{\check{\alpha} \in \tau_f}_\B \in H$ is equivalent to $\alpha \in (\tau_f)^H$ and the latter is equivalent to $b_{C_{\alpha}} \in H$ by the definition of $\tau_f$. Contradiction!

%This ends the proof of (\ref{enum-lem:cont-name3}) of Lemma~\ref{cont-name}.

\item[\ref{enum-lem:cont-name4}] Let $G$ be any $\B$-generic filter over $V$. We will show that for all $\alpha < \kappa$, $\alpha \in \bigl(\tau_{(f_{\tau})}\bigr)^G$ if and only if $\alpha \in \tau^G$. In fact, we have the following equivalences:
\begin{align*}
\alpha \in \bigl(\tau_{(f_{\tau})}\bigr)^G \iff & b_{C_{\alpha}} \in G\\
%				\iff   & G \in C_{\alpha}^{V[G]}\\
%				\iff   & f^{V[G]}_{\tau}(G) (\alpha) = 1\\
				\iff   & \Qp{\check{\alpha} \in \tau }_\B \in G\\
				\iff	& \alpha \in \tau^G,
\end{align*}
where $C_{\alpha} = \{ H \in \text{St}(\B) \mid f_{\tau}(H) (\alpha) =1\}$. 

The first equivalence above follows from the definition of $\tau_{(f_{\tau})}$. The third equivalence above follows from the genericity of $G$.

To see the second equivalence above, we show that $b_{C_{\alpha}} = \Qp{\check{\alpha} \in \tau }_\B$. For this equality, it is enough to see that for all $H\in \text{St}(\B)$, $b_{C_{\alpha}} \in H $ if and only if $\Qp{\check{\alpha} \in \tau }_\B \in H$. But $b_{C_{\alpha}} \in H $ is equivalent to $f_{\tau}(H) (\alpha) =1$ by the definition of $C_{\alpha}$ and the property of $b_{C_{\alpha}}$, and $f_{\tau}(H) (\alpha) =1$ is equivalent to $\Qp{\check{\alpha} \in \tau }_\B \in H$ by the definition of $f_{\tau}$, as desired.
\end{description}
This finishes the proof of Lemma~\ref{cont-name}.
\end{proof}

In Lemma~\ref{cont-name}, working with the product topology on $2^{\kappa}$ is essential. In fact, if we work with the initial segment topology (or bounded topology) on $2^{\kappa}$, then Lemma~\ref{cont-name}~(\ref{enum-lem:cont-name4}) {\it fails}. The basic open sets of $2^{\kappa}$ with the initial segment topology are of the form $U_{s} = \{ x \in 2^{\kappa} \mid x \supseteq s\}$ where $s$ is in $2^{<\kappa}$, most likely {\it not} a finite partial function from $\kappa$ to $2$, in particular it is much harder for an $f\colon \St(\B)\to 2^\kappa$ to be continuous with respect this finer topology on $2^\kappa$: %{\it not} necessarily a finite partial function from $\kappa$ to $2$.
\begin{Ex}
Let $\kappa = \omega_1$ and let $\B$ be the completion of the partial order adding $\omega_1$-many Cohen reals, i.e., $\text{Fn}(\omega_1, 2)$ ordered by reverse inclusion. Let $\tau$ be a $\B$-name for a subset of $\omega_1$ coding the generic $\omega_1$-sequence of Cohen reals. %Then for any continuous function $f$ from $\text{St}(\B)$ to $2^{\kappa}$ with the {\it initial segment topology}, $\Vdash_\B \tau_f \neq \tau$. 
Then for any continuous function $f$ from $\text{St}(\B)$ to $2^{\kappa}$ with the {\it initial segment topology}, $\Vdash_\B \tau_f \neq \tau$.

This occurs since for any such $f$, we can split $\B$ in a maximal antichain $\mathcal{A}=\{ b_s \mid s\in 2^\omega\cap V \}$ in $V$  given by the clopen sets in $\St(\B)$ which are the preimage by $f$ of the clopen set (in the initial segment topology) given by an $s\in 2^\omega$ seen as an initial segment of $2^\kappa$, i.e., $O_{b_s} = f^{-1} (U_s)$.
For any $\B$-generic filter $G$ over $V$, $G\cap\mathcal{A}$ is non-empty, giving that $(\tau_f)^G\restriction\omega=s$ %=f(G)\restriction \omega=s$ 
for some $s\in 2^\omega\cap V$. On the other hand, the first $\omega$-many bits of $\tau$ would define a Cohen real over $V$. 
In particular $(\tau_f)^G\neq\tau^G$ for any such $G$.
\end{Ex}

\begin{comment}
Assume we work with the initial segment topology (or bounded topology) on $2^{\kappa}$, then 
Lemma~\ref{cont-name}(\ref{enum-lem:cont-name4}) {\it fails}: the basic open sets of $2^{\kappa}$ 
with the initial segment topology are of the form 
$U_{s} = \{ x \in 2^{\kappa} \mid x \supseteq s\}$ where $s$ in $2^{<\kappa}$ is most likely
{\it not} a finite partial function from $\kappa$ to $2$, in particular it is much harder for an
$f:\St(\B)\to 2^\kappa$ to be continuous with respect this finer topology on $2^\kappa$:
\begin{example}
Let $\kappa = \omega_1$ and let $\B$ be the completion of the partial order adding $\omega_1$-many Cohen reals. 
Let $\tau$ be a $\B$-name for a subset of $\omega_1$ coding the generic $\omega_1$-sequence of Cohen reals. 
Then for any continuous function $f$ from $\text{St}(\B)$ to $2^{\kappa}$ with 
the {\it initial segment topology}, $\Vdash_\B \tau_f \neq \tau$.
This occurs since for any such $f$ we can split $\B$ in a maximal antichain
$\mathcal{A}=\bp{A_s:s\in 2^\omega\cap V}$ in $V$ 
given by the clopen sets in $\St(\B)$ which are the preimage by $f$ of the clopen set (in the bounded topology) 
given by an 
$s\in 2^\omega$ seen as an initial 
segment of $2^\kappa$.
For any $G$ $V$-generic for $\B$ $G\cap\mathcal{A}$ is non-empty, giving that
$(\tau_f)^G\restriction\omega=f(G)\restriction \omega=s$ 
for some $s\in 2^\omega\cap V$. On the other hand the first $\omega$-many bits of $\tau$ 
would define a Cohen real over $V$. 
In particular $(\tau_f)^G\neq\tau^G$ for any such $G$.
\end{example}
\end{comment}

If $\kappa$ is a regular cardinal and $\B$ is $<$$\kappa$-closed as a partial order, then one obtains Lemma~\ref{cont-name} even with the initial segment topology on $2^{\kappa}$. This was observed by Philipp Schlicht.

 For a $\B$-name $\tau$ for a subset of $\kappa$ and a filter $g$ on $\B$, let $\text{val}(\tau, g)$ be the interpretation of $\tau$ via $g$, i.e., $\text{val}(\tau, g) = \{ \alpha < \kappa \mid [\check{\alpha} \in \tau]_B \in g\}$. Note that if $g$ is a $\B$-generic filter over $V$, then $\text{val}(\tau, g)$ is the same as $\tau^g$, i.e., the standard interpretation of a $\B$-name via a filter.

We further connect $\B$-names for a subset of $\kappa$ and continuous functions from $\text{St}(\B)$ to $2^{\kappa}$ with filters on $\B$ which are generic over (sufficiently) elementary substructures of $V$. Let $Z$ be a set and $g$ be a filter on $\B$. We say $g$ is {\it $(Z,\B)$-generic} if for any dense subset $D$ of $\B$ which is in $Z$, $g \cap D \cap Z \neq \emptyset$. 

\begin{Rem}\label{generic collapse}
Let $Z$ be a $\Sigma_{2024}$-elementary substructure of $V$ such that $\kappa , \B \in Z$ and $|Z| = \kappa \subseteq Z$. Let $g$ be a filter on $\B$ which is $(Z,\B)$-generic. Let $\pi \colon Z \to M_Z$ be the transitive collapse and $\bar{g} = \pi \lq\lq (g \cap Z)$. Then
\begin{enumerate}
\item $\bar{g}$ is a $\pi (\B)$-generic filter over $M_Z$, and

\item for all $\B$-names $\tau$ for a subset of $\kappa$ which are in $Z$, $\text{val} (\tau ,g) = \bar{\tau}^{\bar{g}}$, where $\bar{\tau} = \pi (\tau)$.
\end{enumerate}
\end{Rem}

\begin{Lem}\label{elementary substructure name cont}
Let $Z$ be a $\Sigma_{2024}$-elementary substructure of $V$ such that $\kappa , \B \in Z$ and $|Z| = \kappa \subseteq Z$. Let $g$ be a filter on $\B$ which is $(Z,\B)$-generic. Then
\begin{enumerate}[(1)]
\item if $\tau$ is a $\B$-name for a subset of $\kappa$ with $\tau \in Z$, then for any ultrafilter $G$ on $\B$ with $g \subseteq G$, the function $f_{\tau}(G)$ is the characteristic function of the subset $\text{val}(\tau , g)$ of $\kappa$, and\label{lem:elementary substructure name cont item1} 

\item if $f$ is a continuous function from $\text{St}(\B)$ to $2^{\kappa}$ with $f \in Z$, then for any ultrafilter $G$ on $\B$ with $g \subseteq G$, $f(G)$ is the characteristic function of $\text{val}(\tau_f , g)$.\label{lem:elementary substructure name cont item2}
\end{enumerate}
\end{Lem}

\begin{proof}

${}$

\begin{description}
\item[\ref{lem:elementary substructure name cont item1}]
We show that for all $\alpha < \kappa$, $f_{\tau}(G) (\alpha) = 1$ if and only if $\alpha \in \text{val}(\tau ,g)$. In fact, we have the following equivalences:
\begin{align*}
f_{\tau}(G) (\alpha) = 1 \iff & \Qp{\check{\alpha} \in \tau }_\B \in G\\
			\iff  & \Qp{\check{\alpha} \in \tau }_\B \in g\\
			\iff  & \alpha \in \text{val}(\tau , g).
\end{align*}

The first equivalence above follows from the definition of $f_{\tau}$ in Definition~\ref{translations name cont}. The second equivalence above follows from the fact that $\Qp{\check{\alpha} \in \tau }_\B \in Z$ (since $\kappa , \B , \tau \in Z$ and $\kappa \subseteq Z$) and the $(Z,\B)$-genericity of $g$. The third equivalence above is immediate by the definition of $\text{val}(\tau, g)$. 

%This ends the proof of the first item 1. of Lemma~\ref{elementary substructure name cont}.

\item[\ref{lem:elementary substructure name cont item2}]
We show that for all $\alpha < \kappa$, $\alpha \in \text{val}(\tau_f ,g )$ if and only if $f(G) (\alpha) =1 $. In fact, we have the following equivalences:
\begin{align*}
\alpha \in \text{val}(\tau_f , g) & \iff \Qp{\check{\alpha} \in \tau_f}_\B \in g\\
				  & \iff b_{C_{\alpha}} \in g \\
				  & \iff b_{C_{\alpha}} \in G \\
				  & \iff f (G) (\alpha) =1,
\end{align*}
where $C_{\alpha} = \{ H \in \text{St}(\B) \mid f (H) (\alpha) =1 \}$ and $b_{C_{\alpha}}$ is the unique element of $\B$ such that $O_{(b_{C_{\alpha}})} = C_{\alpha}$. 

The first equivalence above follows from the definition of $\text{val}(\tau_f , g)$. The third equivalence above follows from the fact that $b_{C_{\alpha}} \in Z$ (since $\kappa, \B , f \in Z$ and $\kappa \subseteq Z$) and the $(Z,\B)$-genericity of $g$. The fourth equivalence above follows from the definition of $C_{\alpha}$. 

To see the second equivalence above, we show that $\Qp{\check{\alpha} \in \tau_f}_\B = b_{C_{\alpha}}$. Suppose that they are different. Then one can take a $\B$-generic filter $H$ over $V$ such that $\Qp{\check{\alpha} \in \tau_f}_\B$ is in $H$, but $b_{C_{\alpha}}$ is {\it not} in $H$. However, by the genericity of $H$, $[\check{\alpha} \in \tau_f]_B \in H$ is equivalent to $\alpha \in (\tau_f)^H$ and the latter is equivalent to $b_{C_{\alpha}} \in H$ by the definition of $\tau_f$. Contradiction!

\end{description}
This finishes the proof of Lemma~\ref{elementary substructure name cont}.
\end{proof}

\section{Universally Baire sets in $2^{\kappa}$}

Let $\Gamma$ be a definable class of complete Boolean algebras.
We introduce the definition and the basic theory of universally Baire sets in $2^{\kappa}$ with respect to the class $\Gamma$ (the $\text{uB}^{\Gamma}_{\kappa}$ sets). Our definitions generalize the classical notion of universally Baire sets of reals due to Feng, Magidor, and Woodin~\cite{MR1233821}. 
%In this section, we discuss the definition and the basic theory of universally Baire sets in $2^{\kappa}$ with respect to a class of complete Boolean algebras $\Gamma$ ($\mathsf{uB}^{\Gamma}_{\kappa}$ sets), which generalizes the classical notion of universally Baireness for sets of reals due to Feng, Magidor, and Woodin~\cite{MR1233821}. 
We start with the notions of $\kappa$-meagerness and $\kappa$-Baire property for subsets of a topological space:
\begin{Def}
Let $X$ be a topological space and $A$ be a subset of $X$. 
\begin{enumerate}
\item We say $A$ is {\it nowhere dense} if the closure of the interior of $A$ is empty.

\item We say $A$ is {\it $\kappa$-meager} if $A$ is the union of $\kappa$-many nowhere dense subsets of $X$.

\item We say $A$ has the {\it $\kappa$-Baire property in $X$} if there is an open subset $U$ of $X$ such that the symmetric difference between $A$ and $U$ (denoted by $A \triangle U$) is $\kappa$-meager in $X$.
\end{enumerate}
\end{Def}

Note that if the whole space $X$ is $\kappa$-meager in $X$, then then every subset of $X$ is $\kappa$-meager in $X$ and so it has the $\kappa$-Baire property in $X$. Hence the notion of $\kappa$-Baire property becomes non-trivial only when $X$ is {\it not} $\kappa$-meager in $X$, which is equivalent to the forcing axiom for $\kappa$-many dense subsets of a complete Boolean algebra $\B$ when $X$ is the Stone space of $\B$.
\begin{Def}
Let $P$ be a partial order. The {\it forcing axiom $\FA_{\kappa}(P)$} states that for any $\kappa$-many dense subsets $(D_{\alpha} \mid \alpha < \kappa)$ of $P$, there exists a filter $g$ on $P$ such that $g \cap D_{\alpha} \neq \emptyset$ for all $\alpha < \kappa$. 
For a complete Boolean algebra $\B$, we denote by $\FA_{\kappa} (\B)$ the statement $\FA_{\kappa} (\B^+)$, where $\B^+$ is the partial order $\B \setminus \{ 0_{\B} \}$ induced by the Boolean algebra $\B$. 
\end{Def}

%By $\FA_{\kappa} (\B)$, we mean the following forcing axiom: For any $\kappa$-many dense subsets $(D_{\alpha} \mid \alpha < \kappa )$ of $\B$, there exists a filter $g$ on $\B$ such that $g \cap D_{\alpha} \neq \emptyset$ for all $\alpha < \kappa$.

From now on, when we say $D$ is a dense subset of $\B$, we mean that $D$ is a dense subset of the partial order $\B^+$.
\begin{Rem}\label{meager-FA}
Let $X$ be the Stone space of $\B$, i.e., $X=\text{St}(\B)$. Then
\begin{enumerate}
\item a subset $D$ of $\B$ is dense if and only if the set $U_{D} = \bigcup \{O_b \mid b \in D\}$ is open dense in $X$, where $O_b = \{ G \in X \mid b \in G\}$, and

\item the space $X$ is {\it not} $\kappa$-meager if and only if the axiom $\FA_{\kappa} (\B)$ holds.
\end{enumerate}
\end{Rem}

\begin{comment}
Let us recall the following useful fact on forcing axioms:
\begin{Fact}
Let $\dot{\C}$ be a $\B$-name such that $\Vdash_\B \lq\lq \dot{\C} \text{ is a complete Boolean algebra}"$.
\begin{enumerate}
\item If the axiom $\FA_{\kappa} (\B \ast \dot{\C})$ holds, then both $\FA_{\kappa}(\B)$ and $\Vdash_\B \lq\lq \FA(\check{\kappa} , \dot{\C})"$ hold. 

\item If the axiom $\FA_{\kappa} \bigl(\B \ast \dot{\C} \upharpoonright (b, \dot{c}) \bigr)$ holds for all $(b, \dot{c} ) \in \B \ast \dot{\C}$, then $\FA_{\kappa}(\B \upharpoonright b)$ holds for all $b$ in $\B$ and $\Vdash_\B \lq\lq (\forall c \in \dot{\C} ) \ \FA_{\check{\kappa}} (\dot{\C} \upharpoonright c)"$. 
\end{enumerate}
\end{Fact}
\end{comment}

The next lemma states that a sufficiently generic ultrafilter on $\B$ does {\it not} belong to a $\kappa$-meager subset of $\text{St}(\B)$:
\begin{Lem}\label{meager generic}
Let $Z$ be a $\Sigma_{2024}$-elementary substructure of $V$ such that $\kappa , \B \in Z$ and $|Z| = \kappa \subseteq Z$. Let $g$ be a filter on $\B$ which is $(Z,\B)$-generic. Then if $A$ is a $\kappa$-meager subset of $\text{St}(\B)$ with $A \in Z$, then for all ultrafilters $G$ on $\B$ with $G \supseteq g$, $G$ is {\it not} in $A$.
\end{Lem}

\begin{proof}
Let $Z, g, A, G$ be as in the statement of this lemma. We show that $G$ does not belong to $A$. Since $A$ is $\kappa$-meager and $A$ is in $Z$, there is a sequence $(D_{\alpha} \mid \alpha < \kappa)$ of dense subsets of $\B$ in $Z$ such that $A$ is disjoint from $\bigcap_{\alpha < \kappa} U_{D_{\alpha}}$, where $U_{D_{\alpha}} = \bigcup_{b \in D_{\alpha}} O_b$ and $O_b = \{ H \in \text{St}(\B) \mid b \in H\}$. 

Since $\kappa \subseteq Z$, $D_{\alpha}$ is in $Z$ for all $\alpha < \kappa$. Hence by the $(Z,\B)$-genericity of $g$, $g \cap D_{\alpha} \neq \emptyset$ for all $\alpha < \kappa$. Since $G \supseteq g$, $G$ is in $\bigcap_{\alpha < \kappa} U_{D_{\alpha}}$. Since $A$ is disjoint from $\bigcap_{\alpha < \kappa} U_{D_{\alpha}}$, $G$ is not in $A$, as desired.
\end{proof}

We next introduce the notion of $\B$-Baireness in $2^{\kappa}$:
\begin{Def}
Let $A$ be a subset of $2^{\kappa}$. We say $A$ is {\it $\B$-Baire in $2^{\kappa}$} if for any continuous function $f$ from $\text{St}(\B)$ to $2^{\kappa}$, the set $f^{-1}(A)$ has the $\kappa$-Baire property in $\text{St}(\B)$. 
\end{Def}

By Remark~\ref{meager-FA}, the notion of $\B$-Baireness becomes non-trivial only when the forcing axiom $\FA_{\kappa}(\B)$ holds. So assuming $\FA_{\kappa}(\B)$ is harmless in this context. 

We will characterize the notion of $\B$-Baireness in two ways: One is in terms of a $\B$-name for a subset of $2^{\kappa}$ with some special property. The other is in terms of the notion of trees for closed subsets of $Y^{\kappa}$ for some set $Y$ with some absoluteness property. 

Before going into the characterization of $\B$-Baireness, let us fix some notations and introduce some basic notions.
We say a $\B$-name is {\it a $\B$-name for a subset of $2^{\kappa}$} if it is a collection of ordered pairs $(\tau, b)$ where $\tau$ is a $\B$-name for a subset of $\kappa$ and $b$ is an element of $\B$. Let $\dot{A}$ be a $\B$-name for a subset of $2^{\kappa}$, $Z$ be a set, and $g$ be $(Z, \B)$-generic. We use $\dot{A}_Z^g$ to denote the set $\{ \text{val}(\tau , g) \mid \text{ for some } b\in g$, $(\tau,b) \in Z \cap \dot{A} \}$, where $\text{val} (\tau , g) = \{ \alpha < \kappa \mid \Qp{\check{\alpha} \in \tau }_\B \in g \}$.

Recall that for a set $Y$, we write $\text{Fn}(\kappa , Y)$ for the collection of finite partial functions from $\kappa$ to $Y$. For $t \in \text{Fn}(\kappa , Y)$, let $[t] = \{ x \in Y^{\kappa } \mid x \supseteq t\}$. We topologize the set $Y^{\kappa}$ by the product topology where $Y$ is considered as a discrete space. Hence the family $\{ [t] \mid t \in \text{Fn}(\kappa , Y)\}$ is a basis for this topology.

\begin{Def}
Let $Y$ be a set. We say a subset $T$ of $\text{Fn}(\kappa , Y)$ is a {\it $\text{tree}^{\kappa}$ on $Y$} if $T$ is closed under initial segments, i.e., if $t$ is in $T$ and $s$ is a subsequence of $t$, then $s$ is in $T$ as well. 
\end{Def}

The notion of $\text{tree}^{\kappa}$ generalizes the notion of tree in descriptive set theory and it corresponds to a closed subset of $Y^{\kappa}$ in an expected way: 
\begin{Fact}
Let $T$ be a $\text{tree}^{\kappa}$ $T$ on $Y$. Let $[T] = \{ x \in Y^{\kappa} \mid (\forall u \in [Y]^{<\omega} ) \ x \upharpoonright u \in T\}$. Then the set $[T]$ is a closed subset of $Y^{\kappa}$. Conversely, if $C$ is a closed subset of $Y^{\kappa}$, then the set $T = \{ t \in \text{Fn}(\kappa , Y) \mid [t] \cap C \neq \emptyset\}$ is a $\text{tree}^{\kappa}$ on $Y$ and $[T] = C$.
\end{Fact}

If $Y_1$ and $Y_2$ are sets, then we will identify $(Y_1 \times Y_2)^{\kappa}$ with $Y_1^{\kappa} \times Y_2^{\kappa}$ in a natural way, for example, if $T$ is a $\text{tree}^{\kappa}$ on $Y_1 \times Y_2$, then we regard $[T]$ as a subset of $Y_1^{\kappa} \times Y_2^{\kappa}$ and we use $\text{p} [T]$ for the image of $[T]$ via the projection map from $Y_1^{\kappa} \times Y_2^{\kappa}$ to $Y_1^{\kappa}$. Also if $T$ is a $\text{tree}^{\kappa}$ on $Y_1 \times Y_2$, then we will identify $T$ as a subset of $\text{Fn}(\kappa , Y_1) \times \text{Fn} (\kappa , Y_2)$ in a natural way.

Let $T$ be a $\text{tree}^{\kappa}$ on $Y \times W_1$ and $U$ be a $\text{tree}^{\kappa}$ on $Y \times W_2$. We define the $\text{tree}^{\kappa}$ $T\wedge U$ on $Y \times W_1 \times W_2$ as follows: $(s,t,u) \in (T \wedge U)$ if both $(s,t) \in T$ and $(s,u) \in U$ hold. Then $\text{p}[T \wedge U] = \text{p}[T] \cap \text{p}[U]$ holds. 

For $b$ in $\B$, we set $\B \upharpoonright b = \{ b' \in \B \mid b' \le b\}$. The partial order $\B \upharpoonright b$ can been seen as a complete Boolean algebra by setting $1_{(\B \upharpoonright b)} = b$. 

The next lemma states that for a $\text{tree}^{\kappa}$ $T$, the statement \lq\lq $[T] \neq \emptyset$" is absolute between $V$ and its generic extension via $\B$ if the forcing axiom $\FA_{\kappa}(\B\upharpoonright b)$ holds for all $b$ in $\B$:
\begin{Lem}\label{tree branch absoluteness}
Suppose that for all $b$ in $\B$, the forcing axiom $\FA_{\kappa}(\B \upharpoonright b)$ holds. Let $T$ be a $\text{tree}^{\kappa}$ and $G$ be a $\B$-generic filter over $V$. Then $[T] \neq \emptyset$ in $V$ if and only if $V[G] \vDash \lq\lq [T] \neq \emptyset "$.
\end{Lem}

\begin{proof}
The forward direction is easy: If $x$ is in $[T]$ in $V$, then for all $u \in [\kappa]^{<\omega}$, $x \upharpoonright u $ is in $T$, which stays true in $V[G]$. Hence $x$ is in $[T]$ also in $V[G]$.

For the reverse direction, let us assume that $[T] \neq \emptyset$ in $V[G]$. We will show that $[T]$ is non-empty also in $V$. Since $[T] \neq \emptyset$ in $V[G]$, there are a $b_0$ in $G$ and a $\B$-name $\dot{x}$ such that $b_0 \Vdash_\B \lq\lq \dot{x} \in [\check{T}]"$. For each $\alpha < \kappa$, let $D_{\alpha} = \{ b \in (\B \upharpoonright b_0) \mid \text{$b$ decides the values of $\dot{x}(\check{\alpha})$} \}$. Then $D_{\alpha} $ is dense in $\B \upharpoonright b_0$ for all $\alpha < \kappa$. By the forcing axiom $\FA_{\kappa} (\B \upharpoonright b_0)$, there is a filter $g$ on $\B$ meeting $D_{\alpha}$ for all $\alpha < \kappa$. Let $x$ be the interpretation of $\dot{x}$ via $g$. Then it is easy to see that $x \in [T]$ in $V$ and hence $[T] \neq \emptyset$ in $V$ as well, as desired.
\end{proof}

We are now ready to state the characterization of $\B$-Baireness. Recall that the notion of $\B$-Baireness becomes non-trivial only when the forcing axiom $\FA_{\kappa}(\B)$ holds. Recall also that if $\dot{A}$ is a $\B$-name for a subset of $2^{\kappa}$, $Z$ is a sufficiently elementary substructure of $V$ with $\kappa , \B \in Z$ and $|Z| = \kappa \subseteq Z$, and if $g$ is a $(Z,\B)$-generic filter on $\B$, then $\dot{A}_Z^g = \{ \text{val}(\tau , g) \mid \text{ for some } b\in g, (\tau,b) \in  \dot{A} \cap Z \}$, and let $\pi \colon Z \to M_Z$ be the collapsing map and $\bar{g} = \pi \lq\lq (g \cap Z)$.
\begin{Thm}\label{B-Baire char}
Suppose that the axiom $\FA_{\kappa}(\B \upharpoonright b)$ holds for all $b$ in $\B$ and let $A$ be a subset of $2^{\kappa}$. Then the following are equivalent:
\begin{enumerate}[(a)]
\item \label{enum-thm:B-Baire-char1} $A$ is $\B$-Baire.

\item \label{enum-thm:B-Baire-char2} There is a $\B$-name $\dot{A}$ for a subset of $2^{\kappa}$ such that for all $Z \prec_{\Sigma_{2024}} V$ with $\kappa , \B , A , \dot{A} \in Z$ and $|Z | = \kappa \subseteq Z$ and all $(Z, \B)$-generic filters $g$ on $\B$, the equality $\dot{A}^g_Z = A \cap M_Z[\bar{g}]$ holds.

\item \label{enum-thm:B-Baire-char3} For some set $Y$, there are $\text{tree}^{\kappa}$s $T, U$ on $2 \times Y$ such that $A = \text{p} [T]$ and $\Vdash_\B \lq\lq \text{p}[\check{T}] = 2^{\kappa} \setminus \text{p}[\check{U}]"$.
\end{enumerate}
\end{Thm}

\begin{proof}
We will show the directions from \ref{enum-thm:B-Baire-char1} to \ref{enum-thm:B-Baire-char2}, from \ref{enum-thm:B-Baire-char2} to \ref{enum-thm:B-Baire-char3}, and from \ref{enum-thm:B-Baire-char3} to \ref{enum-thm:B-Baire-char1}.

\begin{description}
\item[\ref{enum-thm:B-Baire-char1}$\Rightarrow$\ref{enum-thm:B-Baire-char2}] Let $\dot{A}$ be the following $\B$-name:
\begin{align*}
\dot{A} = \{ (\tau , b) \mid \text{$\tau$ is a $\B$-name for a subset of $\kappa$ and $O_b \setminus f_{\tau}^{-1}(A)$ is $\kappa$-meager}\},
\end{align*}
where $f_{\tau}$ is as in Definition~\ref{translations name cont}. It is easy to see that $\dot{A}$ is a $\B$-name for a subset of $2^{\kappa}$. We will show that $\dot{A}$ witnesses the second item 2.. Let $Z$ be any $\Sigma_{2024}$-elementary substructure of $V$ with $\kappa , \B , A , \dot{A} \in Z$ and $|Z | = \kappa \subseteq Z$, and $g$ be a filter on $\B$ which is $(Z,\B)$-generic. We will show the following:
\begin{align}
\dot{A}^g_Z = A \cap M_Z[\bar{g}], \tag{$\ast$}\label{equation1}
\end{align}
where $\pi \colon Z \to M_Z$ is the transitive collapse and $\bar{g} = \pi \lq\lq (g \cap Z)$.
\begin{proof}[Proof of equation~{\rm (\ref{equation1})}]
We prove both inclusions as follows:
\begin{description}
\item[$\subseteq$] %We first see that $\dot{A}^g_Z \subseteq A \cap M_Z[\bar{g}]$. 
Let $x$ be any element of $\dot{A}^g_Z$. We show that $x$ is also in $A \cap M_Z [\bar{g}]$. Since $x$ is in $\dot{A}^g_Z$, there is a $(\tau ,b) \in \dot{A} \cap Z$ such that $b $ is in $g$ and $x = \text{val}(\tau ,g )$. First note that $x$ is in $M_Z [\bar{g}]$ because $x = \text{val}(\tau ,g ) = \bar{\tau}^{\bar{g}}$ by Remark~\ref{generic collapse} where $\bar{\tau} = \pi (\tau)$ which is in $M_Z$. Therefore, $x$ is in $M_Z [\bar{g}]$.

We show that $x$ is also in $A$. Since $(\tau ,b )$ is in $\dot{A}$, $\tau$ is a $\B$-name for a subset of $\kappa$ and $O_b \setminus f^{-1}_{\tau} (A)$ is $\kappa$-meager in $\text{St}(\B)$. Since $\tau , b , A \in Z$, the set $O_b \setminus f^{-1}_{\tau} (A)$ is also in $Z$. Hence by Lemma~\ref{meager generic}, for all ultrafilters $G$ on $\B$ with $g \subseteq G$, $G$ does not belong to $O_b \setminus f^{-1}_{\tau} (A)$. Take any such an ultrafilter $G$. Since $b \in g$, $G$ belongs to $f^{-1}_{\tau} (A)$ for all $G\in \text{St}(\B)$ with $g \subseteq G$. Therefore, $f_{\tau} (G ) \in A$. But by Lemma~\ref{elementary substructure name cont}, $f_{\tau} (G ) = \text{val} (\tau , g)$. Hence $x = \text{val} (\tau ,g) \in A$, as was to be shown. %This ends the proof of $\dot{A}^g_Z \subseteq A \cap M_Z[\bar{g}]$

\item[$\supseteq$]  %We next see that $A \cap M_Z[\bar{g}] \subseteq \dot{A}^g_Z$. 
Let $x$ be any element of $A \cap M_Z [\bar{g}]$. We will prove that $x$ is also in $\dot{A}^g_Z$. Since $x$ is in $M_Z [\bar{g}]$, there is a $\pi (\B)$-name $\bar{\tau}$ in $M_Z$ such that $x = \bar{\tau}^{\bar{g}}$ and $\bar{\tau}$ is a $\pi (\B)$-name for a subset of $\kappa$ in $M_Z$ (note that $\pi (\kappa) = \kappa$ because $\kappa \subseteq Z$). Let $\tau$ be in $Z$ such that $\pi (\tau) = \bar{\tau}$. Then $\tau$ is a $\B$-name for a subset of $\kappa$ in $Z$ and hence in $V$ as well. 
\begin{Claim}\label{D dense}
Let $D$ be the following subset of $\B$:
\begin{align*}
D = \{ b \in \B \mid  \text{either $O_b \setminus f^{-1}_{\tau} (A)$ is $\kappa$-meager or $O_b \cap f^{-1}_{\tau}(A)$ is $\kappa$-meager}\}.
\end{align*}
Then $D$ is dense in $\B$.
\end{Claim}

\begin{proof}[Proof of Claim~\ref{D dense}]

Take any element $b$ of $\B$. We will show that there is a $b' \le b $ such that $b' \in D$. Since $A$ is $\B$-Baire, $f^{-1}_{\tau}(A)$ has the $\kappa$-Baire property in $\text{St}(\B)$. Hence there is an open set $U$ in $\text{St}(\B)$ such that the symmetric difference $f^{-}_{\tau} (A) \triangle U$ is $\kappa$-meager in $\text{St}(\B)$. 

Suppose that $O_b \cap U$ is empty. Then it follows that $O_b \cap f^{-1}_{\tau}(A)$ is $\kappa$-meager in $\text{St}(\B)$ and hence $b$ is in $D$, as desired.

Suppose that $O_b \cap U$ is non-empty. Then there is a $b' \le b$ such that $O_{b'} \subseteq (O_b \cap U)$ and it follows that $O_{b'} \setminus f^{-1}_{\tau}(A)$ is $\kappa$-meager and $b'\in D$, as desired.

This finishes the proof of Claim~\ref{D dense}.
\end{proof}

Note that the dense set $D$ in Claim~\ref{D dense} is in $Z$ because $\B, \tau , A$ are all in $Z$. Hence by the $(Z,\B)$-genericity of $g$, there is a $b$ in $D \cap g \cap Z$. Since $b$ is in $D$, either $O_b \setminus f^{-1}_{\tau}(A)$ is $\kappa$-meager or $O_b \cap f^{-1}_{\tau}(A)$ is $\kappa$-meager. 

We argue that $O_b \setminus f^{-1}_{\tau}(A)$ is $\kappa$-meager. Suppose not. Then $O_b \cap f^{-1}_{\tau}(A)$ is $\kappa$-meager. Since $b$ is in $Z$, $O_b \cap f^{-1}_{\tau}(A)$ is also in $Z$. Hence by Lemma~\ref{meager generic}, for all ultrafilters $G$ on $\B$ with $g \subseteq G$, $G$ does not belong to $O_b \cap f^{-1}_{\tau}(A)$. Take such a $G$. Then since $b \in g \subseteq G$, $G$ cannot belong to $f^{-1}_{\tau}(A)$ and hence $f_{\tau}(G) \notin A$. But by Lemma~\ref{elementary substructure name cont} and Remark~\ref{generic collapse}, $f_{\tau}(G) = \text{val}(\tau , g) = \bar{\tau}^{\bar{g}} = x$ and we assumed that $x$ was in $A$. Contradiction! Therefore, $O_b \setminus f^{-1}_{\tau}(A)$ is $\kappa$-meager. 

By the definition of $\dot{A}$, we can conclude that $(\tau , b)$ is in $\dot{A} \cap Z$. Since $b$ is in $g$, by the definition of $\dot{A}^g_Z$, we have $x = \text{val}(\tau , g) \in \dot{A}^g_Z$, as was to be shown. %This finishes the proof of $A \cap M_Z[\bar{g}] \subseteq \dot{A}^g_Z$.
\end{description}
The equation~{\rm (\ref{equation1})} is proved.
\end{proof}

This concludes the proof of \ref{enum-thm:B-Baire-char1}$\Rightarrow$\ref{enum-thm:B-Baire-char2}.

\item[\ref{enum-thm:B-Baire-char2}$\Rightarrow$\ref{enum-thm:B-Baire-char3}] Let us fix a $\B$-name $\dot{A}$ for a subset of $2^{\kappa}$ witnessing \ref{enum-thm:B-Baire-char2}. We will show that there are $\text{tree}^{\kappa}$s witnessing \ref{enum-thm:B-Baire-char3}. The arguments are a combination of Woodin's argument in \cite[Theorem~2.53]{Woodin_2010} and the argument in the paper by Feng, Magidor, and Woodin~\cite{MR1233821}. 

First notice that $\FA_\kappa(\B)$ entails that all cardinals less or equal than $\kappa$ are preserved by forcing with $\B$.

%First notice that $\kappa$ stays to be a cardinal in $V[G]$ because $\FA_{\kappa} (\B)$ holds by assumption. So when we write $|Z|$ for a set $Z$ of size $\kappa$ in $V$, there is no ambiguity on which model ($V$ or $V[G]$) we compute the cardinality of $Z$. 

Second note that if $a$ is any set and if $G$ is a $\B$-generic filter over $V$, then in $V[G]$, there is a $\Sigma_{2024}$-elementary substructure $Z$ of $V$ such that $a \in Z, |Z| = \kappa \subseteq Z$, and $G$ is $(Z,\B)$-generic. In fact, one can take a $Z$ in such a way that $(Z , \in , Z \cap G)$ is $\Sigma_{2024}$-elementary in $(V, \in , G)$ inside $V[G]$, which gives us the condition that $G$ is $(Z, \B)$-generic. 

For each $\B$-name $\tau$ for a subset of $\kappa$, let us fix a $\B$-name $\dot{Z}_{\tau}$ such that $\Vdash_\B \lq\lq \check{\kappa} , \check{B} , \check{A}, \check{\tau},\check{\dot{A}} \in \dot{Z}_{\tau}, \dot{Z}_{\tau} \prec_{\Sigma_{2024}} \check{V}, |\dot{Z}_{\tau}| = \kappa \subseteq \dot{Z}_{\tau}$, and $\dot{G}$ is $(\dot{Z}_{\tau}, \check{B})$-generic", where $\check{V}$ denotes the ground model $V$ in the generic extension, (which by a Theorem of Laver~\cite[Theorem~3]{MR2364192} is a definable class in $V[G]$ in the parameters $\B$, $G$, and $H_{|\B|^+}^V$),  $\dot{G}$ is a canonical $\B$-name for the generic filter, and $\check{\tau}$ is a $\B$-name for $\tau$, i.e., $\check{\tau}^G = \tau$ for any filter $G$ (so the above condition ensures that $\tau$ itself is in $\dot{Z}_{\tau}^G$ and $\tau^G$ may {\it not} be in $\dot{Z}_{\tau}^G$ in general if $G$ is $\B$-generic over $V$). The same applies to $\check{\dot{A}}$. 
Let $\dot{h}_{\tau}$ be a $\B$-name such that $\Vdash_\B \lq\lq \dot{h}_{\tau}$ is a surjection from $\check{\kappa}$ to $\dot{Z}_{\tau}$". Let $R_{\tau}$ be the relation coding $\dot{h}_{\tau}$, i.e., $R_{\tau} = \{ (b , \alpha , z ) \mid b \in \B, \alpha < \kappa , z \in V, \text{ and } b \Vdash_B \lq\lq \dot{h}_{\tau} (\check{\alpha}) = \check{z}" \}$. 

Let $\gamma$ be a sufficiently large ordinal such that $V_{\gamma} \prec_{\Sigma_{2024}} V$ and $V_{\gamma}$ contains all the sets we have considered so far in the arguments for the direction from \ref{enum-thm:B-Baire-char2} to \ref{enum-thm:B-Baire-char3} in this theorem.

Let $\mathcal{D}_{\tau}$ be the collection of all dense subsets of $\B$ which are definable in the structure $(V_{\gamma}, \in , R_{\tau})$ with parameters from the set $\{ \kappa , \B , A , \tau , \dot{A} \} \cup \kappa$. Note that the size of $\mathcal{D}_{\tau}$ is at most $\kappa$ and hence by the axiom $\FA_{\kappa} (\B)$, there is a filter $g$ on $\B$ meeting all dense sets in $\mathcal{D}_{\tau}$. Let $Z^g_{\tau} = \{ z \in V \mid \text{ for some } b \in g \text{ and } \alpha < \kappa,  (b, \alpha , z) \in R_{\tau}\}$. Then by the definition of $\mathcal{D}_{\tau}$ and the $\mathcal{D}_{\tau}$-genericity of $g$, the set $Z^g_{\tau}$ is a $\Sigma_{2024}$-elementary substructure of $V$ such that $\kappa , \B, A, \dot{A}, \tau \in Z^g_{\tau}$,  $|Z^g_{\tau}| = \kappa \subseteq Z^g_{\tau}$, and that $g$ is $(Z^g_{\tau} ,\B)$-generic. 

We shall define the desired $\text{tree}^{\kappa}$s $T$ and $U$ witnessing the third item 3. using the assignment $\tau \mapsto \mathcal{D}_{\tau} $. For each $\B$-name $\tau$ for a subset of $\kappa$, let $(D^{\tau}_{\alpha} \mid \alpha < \kappa)$ be an enumeration of $\mathcal{D}_{\tau}$. The $\text{tree}^{\kappa}$ $T$ will be the disjoint union $\coprod_{\tau} T_{\tau}$ where $\tau$ ranges over all the $\B$-names for a subset of $\kappa$ and let $T_{\tau}$ be the collection of pairs $(s, t) \in \text{Fn}(\kappa , 2) \times \text{Fn} (\kappa , \B)$ such that
\begin{enumerate}
\item for all $\alpha \in \text{dom}(t)$, $t(\alpha ) \Vdash_\B \lq\lq \tau \in \dot{A}"$,

\item $\text{dom}(s) = \text{dom}(t)$,

\item for all $\alpha \in \text{dom}(s)$, $s(\alpha ) =1 \iff t(\alpha) \Vdash_\B \lq\lq \check{\alpha} \in \tau "$,

\item for all $\alpha \in \text{dom}(s)$, $s(\alpha ) =0 \iff t(\alpha) \Vdash_\B \lq\lq \check{\alpha} \notin \tau "$,

\item for all $\alpha \in \text{dom}(t)$, $t(\alpha) \in D^{\tau}_{\alpha}$, and
\item $\bigwedge \{ t(\alpha) \mid \alpha \in \text{dom}(t)\} \neq 0_B$.
\end{enumerate}

Note that $T_{\tau_1}$ and $T_{\tau_2}$ could be the same as a set even if $\tau_1 \neq \tau_2$, but we will distinguish them when we consider the disjoint union of them for obtaining $T$. It is easy to see that $T_{\tau}$ is a $\text{tree}^{\kappa}$ on $2 \times B$. Hence $T$ is also a $\text{tree}^{\kappa}$.

The definition of $U$ is the same as $T$ except the first item above for the definition of $T_{\tau}$: The $\text{tree}^{\kappa}$ $U$ will be the disjoint union $\coprod_{\tau} U_{\tau}$ where $\tau$ ranges over all the $\B$-names for a subset of $\kappa$ and let $U_{\tau}$ be the collection of pairs $(s, u) \in \text{Fn}(\kappa , 2) \times \text{Fn} (\kappa , \B)$ such that
\begin{enumerate}
\item for all $\alpha \in \text{dom}(u)$, $u(\alpha ) \Vdash_\B \lq\lq \tau \notin \dot{A}"$,

\item $\text{dom}(s) = \text{dom}(u)$,

\item for all $\alpha \in \text{dom}(s)$, $s(\alpha ) =1 \iff u(\alpha) \Vdash_\B \lq\lq \check{\alpha} \in \tau "$,

\item for all $\alpha \in \text{dom}(s)$, $s(\alpha ) =0 \iff u(\alpha) \Vdash_\B \lq\lq \check{\alpha} \notin \tau "$,

\item for all $\alpha \in \text{dom}(u)$, $u(\alpha) \in D^{\tau}_{\alpha}$, and

\item $\bigwedge \{ u(\alpha) \mid \alpha \in \text{dom}(u)\} \neq 0_B$.
\end{enumerate}

We will show that $\text{p}[T] = A$ and $\Vdash_\B \lq\lq \text{p}[\check{T}] = 2^{\kappa} \setminus \text{p} [\check{U}]"$.
This suffices to complete the proof of the direction \ref{enum-thm:B-Baire-char2}$\Rightarrow$\ref{enum-thm:B-Baire-char3}.

We first prove 
\begin{align}
\text{p}[T] = A. \tag{+}\label{equation2}
\end{align}
\begin{proof}
${}$

\begin{description}
\item[$\subseteq$] %We first argue that $\text{p}[T] \subseteq A$. 
Take any element $x$ of $\text{p}[T]$. We will show that $x$ is also in $A$. Since $x$ is in $\text{p}[T]$, there are a $\B$-name $\tau$ for a subset of $\kappa$ and a $\vec{g} \in \B^{\kappa}$ such that $(x,\vec{g}) \in [T_{\tau}]$. By the sixth item of the definition of $T_{\tau}$, it follows that the range of $\vec{g}$ generates a filter. Let $g$ be the filter generated by the range of $\vec{g}$. Then by the fifth item of the definition of $T_{\tau}$, $g$ meets all the dense sets in $\mathcal{D}_{\tau}$ and hence $Z^g_{\tau}$ is a $\Sigma_{2024}$-elementary substructure of $V$ such that $\kappa , \B, A, \dot{A}, \tau \in Z^g_{\tau}$,  $|Z^g_{\tau}| = \kappa \subseteq Z^g_{\tau}$, and that $g$ is $(Z^g_{\tau} ,\B)$-generic, where $Z^g_{\tau} = \{ z \in V \mid \text{ for some } b \in g \text{ and } \alpha < \kappa,  (b, \alpha , z) \in R_{\tau}\}$. For ease of notation, set $Z = Z^g_{\tau}$. By the first item of the definition of $T_{\tau}$, it follows that $\text{val}(\tau , g) \in \dot{A}_Z^g$. By \ref{enum-thm:B-Baire-char2}, $\dot{A}_Z^g = A \cap M_Z [\bar{g}]$ and hence $\text{val}(\tau , g) \in A$. But by the third and fourth items of the definition of $T_{\tau}$, $x = \text{val} (\tau , g)$. Therefore, $x $ is in $A$, as desired. %This finishes the proof of $\text{p}[T] \subseteq A$. 

\item[$\supseteq$] %We next argue that $A \subseteq \text{p}[T]$. 
Take any $x$ in $A$. We will show that $x$ is also in $\text{p}[T]$. Let $\tau = \check{x}$. Then since $\mathcal{D}_{\tau}$ is of cardinality at most $\kappa$, by the forcing axiom $\FA_{\kappa} (\B)$, there is a filter $g$ on $\B$ meeting all the dense sets in $\mathcal{D}_{\tau}$. Hence $Z^g_{\tau}$ is a $\Sigma_{2024}$-elementary substructure of $V$ such that $\kappa , \B, A, \dot{A}, \tau \in Z^g_{\tau}$,  $|Z^g_{\tau}| = \kappa \subseteq Z^g_{\tau}$, and that $g$ is $(Z^g_{\tau} ,\B)$-generic, where $Z^g_{\tau} = \{ z \in V \mid \text{ for some } b \in g \text{ and } \alpha < \kappa,  (b, \alpha , z) \in R_{\tau}\}$. For ease of notation, set $Z = Z^g_{\tau}$. By \ref{enum-thm:B-Baire-char2}, we have $A \cap M_Z[\bar{g}] = \dot{A}_Z^g$. Since $x$ is in $A$ and $x = \check{x}^g = \text{val}(\tau ,g) = \bar{\tau}^{\bar{g}} \in M_Z[\bar{g}]$ (by Remark~\ref{generic collapse}), $x$ is also in $\dot{A}_Z^g$. But $x = \check{x}^g = \text{val}(\tau , g)$. Hence $\text{val}(\tau , g) $ is also in $\dot{A}_Z^g$. By the $(Z,\B)$-genericity and $\mathcal{D}_{\tau}$-genericity of $g$, for each $\alpha < \kappa$, one can take a $b_{\alpha} \in g \cap Z$ such that $b_{\alpha} \Vdash_\B \lq\lq \tau \in \dot{A}"$, $b_{\alpha} \in D^{\tau}_{\alpha}$, and that $ x( \alpha) =1 \iff b_{\alpha} \Vdash_\B \lq\lq \check{\alpha} \in \tau"$. Let $\vec{g} \in \B^{\kappa}$ be such that $\vec{g}(\alpha) = b_{\alpha}$ for each $\alpha < \kappa$. Then it is easy to see that $(x, \vec{g})$ is in $[T]$ and hence $x $ is in $\text{p}[T]$, as desired.
\end{description}
This finishes the proof of $\text{p}[T] = A$.
\end{proof}

We next prove that $\Vdash_\B \lq\lq \text{p}[\check{T}] = 2^{\kappa} \setminus \text{p} [\check{U}]"$. Let $G$ be any $\B$-generic filter over $V$. We will show that in $V[G]$, $\text{p}[T] = 2^{\kappa} \setminus \text{p}[U]$ holds. For this, it is enough to prove that in $V[G]$, $\text{p}[T] \cup \text{p}[U] = 2^{\kappa}$ and $\text{p}[T] \cap \text{p}[U] = \emptyset$. 

We first argue that $\text{p}[T] \cup \text{p}[U] = 2^{\kappa}$ in $V[G]$. Let $x$ be any element of $2^{\kappa}$ in $V[G]$. We will show that $x$ is either in $\text{p}[T]$ or in $\text{p}[U]$. Let $\tau$ be a $\B$-name for a subset of $\kappa$ in $V$ such that $\tau^G = x$. Let us assume that $x$ is in $\dot{A}^G$ . We will prove that $x$ is in $\text{p}[T_{\tau}]$ in $V[G]$ (one can prove that $x$ is in $\text{p}[U_{\tau}]$ in $V[G]$ in the same way in case $x$ is {\it not} in $\dot{A}^G$). Since $x  = \tau^G \in \dot{A}^G$ and $G$ is $\B$-generic over $V$, in $V[G]$, for each $\alpha < \kappa$, one can take a $b_{\alpha} \in G$ such that $b_{\alpha} \Vdash_\B \lq \lq \tau \in \dot{A}"$ in $V$, $b_{\alpha} \in D^{\tau}_{\alpha}$, and that $x(\alpha) =1$ if and only if $b_{\alpha} \Vdash_\B \lq\lq \check{\alpha} \in \tau"$ holds in $V$. In $V[G]$, let $\vec{g} \in \B^{\kappa}$ be such that $\vec{g}(\alpha) = b_{\alpha}$ for each $\alpha < \kappa$. Then it is easy to see that $(x , \vec{g}) \in \text{p}[T_{\tau}]$ in $V[G]$, as desired. This finishes the proof of $\text{p}[T] \cup \text{p}[U] = 2^{\kappa}$ in $V[G]$.

We next argue that $\text{p}[T] \cap \text{p}[U] = \emptyset$ in $V[G]$. Suppose not. Then there is an $x$ which is in both $\text{p}[T]$ and $\text{p}[U]$. We will derive a contradiction. In $V$, let $T \wedge U$ be the $\text{tree}^{\kappa}$ on $2 \times B \times B$ such that $(s,t,u) \in (T \wedge U)$ if both $(s,t) \in T$ and $(s,u) \in U$ hold. Then both in $V$ and $V[G]$, $\text{p}[T \wedge U] = \text{p}[T] \cap \text{p}[U]$ holds. Therefore, in $V[G]$, $\text{p}[T \wedge U]$ and hence $[T\wedge U]$ is non-empty. Since the forcing axiom $\FA_{\kappa} (\B \upharpoonright b)$ holds for all $b$ in $\B$, by Lemma~\ref{tree branch absoluteness}, $[T \wedge U]$ is non-empty also in $V$.

We now derive a contradiction for the assumption that $\text{p}[T] \cap \text{p}[U] \neq \emptyset$ in $V[G]$ by using the fact that $[T \wedge U]$ is not empty also in $V$. We work in $V$. Since $[T \wedge U]$ is nonempty, both $[T]$ and $[U]$ are non-empty. Hence there are $x, \vec{g}, \vec{h}, \sigma$, and $\tau$ such that $(x, \vec{g} ) \in [T_{\sigma}]$ and $(x, \vec{h}) \in [U_{\tau}]$. By the forcing axiom $\FA_{\kappa} (\B)$, there is a filter $g$ on $\B$ meeting all the dense sets in both $\mathcal{D}_{\sigma}$ and $\mathcal{D}_{\tau}$. Let $Z_1 = Z^{\sigma}_g$ and $Z_2 = Z^{\tau}_h$. Then as we argued for $\text{p}[T] \subseteq A$, $x$ is in $\dot{A}_{Z_1}^g$ and {\it not} in $\dot{A}_{Z_2}^g$ while $\dot{A}_{Z_1}^g = A \cap M_{Z_1}[\bar{g}]$ and $\dot{A}_{Z_2}^g = A \cap M_{Z_2}[\bar{g}]$, which implies that $x$ is in $A$ and {\it not} in $A$. Contradiction! Therefore, $\text{p}[T] \cap \text{p}[U] = \emptyset$ in $V[G]$, as desired.

This concludes the proof of \ref{enum-thm:B-Baire-char2}$\Rightarrow$\ref{enum-thm:B-Baire-char3}.

\item[\ref{enum-thm:B-Baire-char3}$\Rightarrow$\ref{enum-thm:B-Baire-char1}]
Let $T, U$ be $\text{tree}^{\kappa}$s witnessing \ref{enum-thm:B-Baire-char3}. Let $f$ be any continuous function from $\text{St}(\B)$ to $2^{\kappa}$. We will show that $f^{-1}(A)$ has the $\kappa$-Baire property in $\text{St}(\B)$. Let $\tau_f$ be as in Definition~\ref{translations name cont} and let $b_0 = \Qp{\tau_f \in \text{p}[\check{T}] }_\B$ and $b_1 = \Qp{\tau_f \in \text{p}[\check{U}] }_\B$. Since $\Vdash_\B \lq\lq \text{p}[\check{T}] = 2^{\kappa} \setminus \text{p}[\check{U}]"$, it follows that $b_0 \vee b_1 = 1_B$ and $b_0 \wedge b_1 = 0_B$. Therefore, it is enough to argue that both $O_{b_0} \setminus f^{-1} (A)$ and $O_{b_1} \cap f^{-1}(A) \bigl( = O_{b_1} \setminus (2^{\kappa} \setminus f^{-1}(A))\bigl)$ are $\kappa$-meager in $\text{St}(\B)$. 

We only show that $O_{b_0} \setminus f^{-1}(A)$ is $\kappa$-meager in $\text{St}(\B)$ (by the same argument, one can show that $O_{b_1} \cap f^{-1}(A)$ is $\kappa$-meager in $\text{St}(\B)$ as well). We will use the same argument as in the direction \ref{enum-thm:B-Baire-char2}$\Rightarrow$\ref{enum-thm:B-Baire-char3}, namely taking a \lq good' $\B$-name for a sufficiently elementary substructure of $V$ (such as $\dot{Z}_{\tau}$) and defining a suitable collection of dense subsets of $\B$ (such as $\mathcal{D}_{\tau}$) from such a $\B$-name. 

Let us fix a $\B$-name $\dot{Z}$ such that $\Vdash_\B \lq\lq \check{\kappa} , \check{B} , \check{A}, \check{f}, \check{T}, \check{U} \in \dot{Z}, \dot{Z} \prec_{\Sigma_{2024}} \check{V}, |\dot{Z}| = \kappa \subseteq \dot{Z}$, and $\dot{G}$ is $(\dot{Z}, \check{B})$-generic". %, where $\check{V}$ denotes the ground model $V$ in a generic extension. 
Let $\dot{h}$ be a $\B$-name such that $\Vdash_\B \lq\lq \dot{h}$ is a surjection from $\check{\kappa}$ to $\dot{Z}$". Let $R$ be the relation coding $\dot{h}$, i.e., $R = \{ (b , \alpha , z ) \mid b \in \B, \alpha < \kappa , z \in V, \text{ and } b \vDash_B \lq\lq \dot{h} (\check{\alpha}) = \check{z}" \}$. 

Let $\gamma$ be a sufficiently large ordinal such that $V_{\gamma} \prec_{\Sigma_{2024}} V$ and $V_{\gamma}$ contains all the sets we have considered so far in the arguments for the direction \ref{enum-thm:B-Baire-char3}$\Rightarrow$\ref{enum-thm:B-Baire-char1}.

Let $\mathcal{D}$ be the collection of all dense subsets of $\B$ which are definable in the structure $(V_{\gamma}, \in , R)$ with parameters from the set $\{ \kappa , \B , A , T, U , f\} \cup \kappa$. 

Note that the size of $\mathcal{D}$ is at most $\kappa$. Hence by Remark~\ref{meager-FA}, it is enough to prove that if $G$ is an ultrafilter on $\B$ meeting all the dense sets in $\mathcal{D}$ such that $b_0$ is in $G$, then $f(G) \in A$. Let $Z^G = \{ z \in V \mid \text{ for some } b \in G \text{ and } \alpha < \kappa,  (b, \alpha , z) \in R\}$. Then by the definition of $\mathcal{D}$ and the $\mathcal{D}$-genericity of $G$, the set $Z^G$ is a $\Sigma_{2024}$-elementary substructure of $V$ such that $\kappa , \B, A, T, U, f \in Z^G$,  $|Z^G| = \kappa \subseteq Z^G$, and that $G$ is $(Z^G ,\B)$-generic. Since $\tau_f$ is simply definable from $f$ in $V$, $\tau_f$ is also in $Z^G$. Since $b_0 \in G$, $b_0 = [ \tau_f \in \text{p}[\check{T}]]_B$, $\kappa \subseteq Z^G$, and $G$ is $(Z^G, B)$-generic, it follows that $\text{val}(\tau_f , G)$ is in $\text{p}[T]$. But by Lemma~\ref{elementary substructure name cont}, $f(G) = \text{val}(\tau_f , G)$. Therefore, $f(G) \in \text{p}[T] = A$, as desired. This finishes the proof of the statement that $O_{b_0} \setminus f^{-1}(A)$ is $\kappa$-meager in $\text{St}(\B)$ and hence of the direction \ref{enum-thm:B-Baire-char3}$\Rightarrow$\ref{enum-thm:B-Baire-char1}.
\end{description}

This finishes the proof of Theorem~\ref{B-Baire char}.
\end{proof}

%The arguments for Theorem~\ref{B-Baire char} also give us that in  of the theorem, 
In Theorem~\ref{B-Baire char}~\ref{enum-thm:B-Baire-char3}, the set $\text{p}[T]$ in a generic extension via $\B$ does not depend on the choice of the pair $(T,U)$:
\begin{Lem}\label{tree invariance}
Suppose that for all $b$ in $\B$, the forcing axiom $\FA_{\kappa}(\B)$ holds. Let $A$ be a $\B$-Baire subset of $2^{\kappa}$ witnessed by two different pairs $(T_1, U_1)$, $(T_2, U_2)$ of $\text{tree}^{\kappa}$s. Then $\Vdash_\B \lq\lq \text{p}[\check{T_1}] = \text{p}[\check{T_2}]"$.
\end{Lem}

\begin{proof}
Let $(T_1, U_1)$, $(T_2, U_2)$ be pairs of $\text{tree}^{\kappa}$s witnessing that $A$ is $\B$-Baire. We will show that $\Vdash_\B \lq\lq \text{p}[T_1] = \text{p}[T_2]"$. 

Recall that if $T$ is a $\text{tree}^{\kappa}$ on $Y \times W_1$ and $U$ is a $\text{tree}^{\kappa}$ on $Y \times W_2$, then $T \wedge U$ is the $\text{tree}^{\kappa}$ on $Y \times W_1 \times W_2$ such that $(s,t,u) \in (T\wedge U)$ if both $(s,t) \in T$ and $(s,u) \in U$ hold. Then $\text{p}[T \wedge U] = \text{p}[T] \cap \text{p}[U]$ holds. 

Since both $(T_1, U_1)$ and $(T_2, U_2)$ witness that $A$ is $\B$-Baire, for each $i =0,1$, $\text{p}[T_i] = A$ and $\text{p}[T_i] = 2^{\kappa} \setminus \text{p}[U_i]$ hold in $V$. In particular, $\text{p} [T_1] \cap \text{p}[U_2] = \emptyset$ and hence $[T_1 \wedge U_2] = \emptyset$ in $V$. Since the forcing axiom $\FA_{\kappa} (\B \upharpoonright b)$ holds for all $b$ in $\B$, by Lemma~\ref{tree branch absoluteness}, $\Vdash_\B \lq\lq [\check{T_1} \wedge \check{U_2}] = \emptyset"$ and hence $\Vdash_\B \lq\lq \text{p}[\check{T_1}] \cap \text{p} [\check{U_2}] = \emptyset"$. 

By the same argument, one can prove that $\Vdash_\B \lq \lq \text{p}[\check{T_2}] \cap \text{p}[\check{U_1}] = \emptyset "$. Since $\Vdash_\B \lq\lq \text{p}[\check{T_1}] = 2^{\kappa} \setminus \text{p}[\check{U_1}] \text{ and } \text{p} [\check{T_2}] = 2^{\kappa} \setminus \text{p}[\check{U_2}]"$, all in all, one can conclude that $\Vdash_\B \lq\lq \text{p}[\check{T_1}] = \text{p}[\check{T_2}]"$, as desired.
\end{proof}

Using Lemma~\ref{tree invariance}, one can naturally interpret a $\B$-Baire subset of $2^{\kappa}$ in a generic extension via $\B$:
\begin{Def}\label{B-Baire lift up}
Suppose that for all $b$ in $\B$, the forcing axiom $\FA_{\kappa}(\B)$ holds. Let $A$ be a $\B$-Baire subset of $2^{\kappa}$ and $G$ be a $\B$-generic filter over $V$. Then, in $V[G]$, let $A^G = \text{p}[T]$, where $T$ is a $\text{tree}^{\kappa}$ in $V$ such that for some $\text{tree}^{\kappa}$ $U$ in $V$, the pair $(T,U)$ witnesses that $A$ is $\B$-Baire in $2^{\kappa}$. We call $A^G$ a {\it natural interpretation of $A$ in $V[G]$}.
\end{Def}

Note that by Lemma~\ref{tree invariance}, the definition of $A^G$ does not depend on the choice of pair $(T,U)$.

Under the presence of Woodin cardinals, one can give another characterization of $\B$-Baire sets in terms of generic absoluteness and stationary tower correctness as in \cite[Theorem~3.3.7]{MR2069032}. For the proof of this characterization, we will use the theory of stationary tower forcing developped by Woodin. The reader can conslut Larson's book~\cite{MR2069032} for the theory of stationary tower forcing. For a Woodin cardinal $\delta$ above $\kappa$, we write $P^{\kappa}_{<\delta}$ for the stationary tower forcing whose conditions $S$ in $V_{\delta}$ are stationary sets in $\bigcup S$ such that $S$ consists of $Z$ such that $|Z| = \kappa \subseteq Z$. If $H$ is $P^{\kappa}_{<\delta}$-generic over $V$, then as in Larson's book~\cite{MR2069032}, $H$ induces a generic embedding $j \colon V \to M \subseteq V[H]$ such that the critical point of $j$ is $\kappa^+$ in $V$, and $M$ is closed under $\kappa$-sequences in $V[H]$.

\begin{Thm}\label{uB generic absoluteness}
Suppose that the axiom $\FA_{\kappa} (\B \upharpoonright b)$ holds for all $b$ in $\B$. Let $\delta$ be a Woodin cardinal above $\kappa$ and $\text{rank}(\B)$, and let $A$ be a subset of $2^{\kappa}$. Then the following are equivalent:
\begin{enumerate}[(a)]
\item \label{uB-generic-absoluteness1} $A$ is $\B$-Baire, and 

\setcounter{enumi}{3}
\item \label{uB-generic-absoluteness2} there are a formula $\phi (v_1, v_2)$ and a set $a$ such that $A = \{ x \in 2^{\kappa} \mid \phi [x,a]\}$ and that for all $\B$-generic filters $G$ over $V$ and all $P^{\kappa}_{<\delta}$-generic filters $H$ over $V$ such that $G \in V[H]$, and all $x \in V[G] \cap 2^{\kappa}$, if $j \colon V \to M$ is the generic embedding induced by $H$, then 
\begin{align*}
V[G] \vDash \phi [x,a] \iff M \vDash \phi [x, j(a)].
\end{align*}
\end{enumerate}
\end{Thm}

\begin{proof}
${}$

\begin{description}
\item[\ref{uB-generic-absoluteness1}$\Rightarrow$\ref{uB-generic-absoluteness2}] 
Let $A$ be $\B$-Baire. We will prove \ref{uB-generic-absoluteness2}. By Theorem~\ref{B-Baire char}, one can take $\text{tree}^{\kappa}$s $T$ and $U$ satisfying Theorem~\ref{B-Baire char}~\ref{enum-thm:B-Baire-char3}. Let $a = T$ and $\phi$ be a formula such that $\phi [x ,T]$ holds if and only if $x$ is in $\text{p}[T]$. We will show that these $\phi$ and $a$ witness \ref{uB-generic-absoluteness2}.

The condition $A = \{ x \mid \phi [x,T]\}$ follows immediately by the property of $T$ such that $A = \text{p}[T]$. We will prove the other condition in \ref{uB-generic-absoluteness2}. Let $G$ be $\B$-generic over $V$ and $H$ be $P^{\kappa}_{<\delta}$-generic over $V$ such that $G \in H$. Let $j \colon V \to M$ be the generic embedding induced by $H$. Note that the critical point of $j$ is $\kappa^{+}$ in $V$ and $M$ is closed under $\kappa$-seuqneces in $V[H]$. We will argue that for all $x$ in $2^{\kappa} \cap V[G]$, $V[G] \vDash \lq\lq x \in \text{p}[T]"$ if and only if $M \vDash \lq\lq x \in \text{p}[j(T)]"$. 

We show the forward direction: If $x $ is in $\text{p}[T]$ in $V[G]$, then there is a $y$ such that $(x, y) \in [T]$. Then since $j(\kappa) = \kappa$, $(x , j \lq\lq y) \in [j(T)]$ in $V[H]$. Since $M$ is closed under $\kappa$-sequences, $(x , j\lq\lq y)$ is also in $M$ and hence $(x,j\lq\lq y) \in \text{p}[j(T)]$ in $M$, therefore $x$ is in $\text{p}[j(T)]$ in $M$, as desired.

We next argue the backward direction by showing its contrapositive, namely if $x$ is {\it not} in $\text{p}[T]$ in $V[G]$, then $x$ is also {\it not} in $\text{p}[j(T)]$ in $M$. Suppose that $x$ is not in $\text{p}[T]$ in $V[G]$. Then by the proerty of $(T,U)$, $x$ is in $\text{p}[U]$ in $V[G]$. Then by the same argument as in the last paragraph, $x$ is also in $\text{p}[j(U)]$ in $M$. However, since $j\colon V \to M$ is elementary and $\text{p}[T] \cap \text{p}[U] = \emptyset$ in $V$, $\text{p}[j(T)] \cap \text{p}[j(U)] = \emptyset $ in $M$ as well. Hence $x $ is not in $\text{p}[j(T)]$ in $M$, as desired.

This completes the proof of the direction \ref{uB-generic-absoluteness1}$\Rightarrow$\ref{uB-generic-absoluteness2}.

%This finishes the proof of the direction from the first item 1. to the second item 2. 

\item[\ref{uB-generic-absoluteness2}$\Rightarrow$\ref{uB-generic-absoluteness1}]  Let $\phi$ and $a$ be as in \ref{uB-generic-absoluteness2}. We will prove that $A$ admits the tree representation as in Theorem~\ref{B-Baire char}~\ref{enum-thm:B-Baire-char3} and hence $A$ is $\B$-Baire by Theorem~\ref{B-Baire char}. 

Let $\dot{A} = \{ (\tau , b) \mid \text{ $\tau$ is a $\B$-name for a subset of $\kappa$}, \Qp{\phi [\tau , \check{a}}_\B = b\}$. Let $\gamma$ be a cardinal such that $\phi$ is absolute between $V_{\gamma}$ and $V$ and $\kappa , \B , A, \dot{A} , a \in V_{\gamma}$. Let $S_0$ be the collection of elementary substructures $Z$ of $V_{\gamma}$ such that $|Z| = \kappa \subseteq Z$ and $\kappa , \B , A , \dot{A} , a \in Z$. Let $S$ be the collection of \lq good' elementary substructures of $V_{\gamma}$ in $S_0$, i.e., $S = \{ Z \in S_0 \mid \text{ for all $(Z, \B)$-generic filters $g$}, \dot{A}^g_Z = A \cap M_Z [\bar{g}]\}$, where $\dot{A}^g_Z , M_Z$, and $\bar{g}$ are as in the paragraph right before Theorem~\ref{B-Baire char}. The next claim states that for club many structures $Z$ in $S_0$, $Z$ is \lq good'.

\begin{Claim}\label{name generic absoluteness}
The set $S_1 = S_0 \setminus S$ is nonstationary in $\bigcup S_0$.
\end{Claim}

\begin{proof}[Proof of Claim~\ref{name generic absoluteness}]
Suppose not. Then $S_1$ is stationary in $\bigcup S_0$. Let $\mu $ be an inaccessible cardinal between $\kappa$ and $\delta$ such that $\B, A$, and $\dot{A}$ are in $V_{\mu}$. Let $\tilde{S_1} = \{ Z \cap V_{\mu} \mid Z \in S_1\}$. Then $\tilde{S_1}$ is stationary in $V_{\mu}$ as well. 

Let $H$ be a $P^{\kappa}_{<\delta}$-generic filter over $V$ such that $\tilde{S_1} $ is in $H$ and let $j\colon V \to M$ be the elementary embedding induced from $H$. Note that the critical point of $j$ is $\kappa^+$ in $V$ and $M$ is closed under $\kappa$-sequences in $V[H]$. 

Since $\tilde{S_1}$ is in $H$, $j\lq\lq V_{\mu}$ is in $j(\tilde{S_1})$. Hence there is a $Z \in j(S_1)$ such that $j\lq\lq V_{\mu} = Z \cap j(V_{\mu})$. Since $Z$ is in $j(S_1)$, by elementarity of $j$, in $M$, there is a $\bigl(Z, j(\B)\bigr)$-generic filter $g$ such that $\dot{A}^g_Z \neq j(A) \cap M_{Z} [\bar{g}]$.  Let $G= \bar{g}$. Then since $\mu $ is large enough for $\B$ and $\dot{A}$, it follows that $\dot{A}_Z^g = \dot{A}_{j\lq\lq V_{\mu}}^g = \{ \text{val}\bigl(j(\tau) , g\bigr) \mid (\exists b ) \ \bigl( j(\tau) , j(b) \bigr) \in j(\dot{A}) \text{ and } j(b) \in g \} = \{ \tau^G \mid (\exists b \in G) \ (\tau , b) \in \dot{A}\} = \dot{A}^G$, while $j(A) \cap M_Z[\bar{g}] = j(A) \cap V_{\mu} [G] = j(A) \cap V[G]$. Therefore, the above inequality implies that $\dot{A}^G \neq j(A) \cap V[G]$.

However, by \ref{uB-generic-absoluteness2}, for all $x \in 2^{\kappa} \cap V[G]$, $V[G] \vDash \lq\lq \phi [x, a]" $ if and only if $M \vDash \lq\lq \phi [x, j(a)]"$, which is equilvalent to $x \in \dot{A}^G \iff x \in j(A)$. Hence $\dot{A}^G = j(A) \cap V[G]$. Contradiction!

This finishes the proof of Claim~\ref{name generic absoluteness}.
\end{proof}

By Claim~\ref{name generic absoluteness}, one can easily see that for any $Z \prec_{\Sigma_{2024}} V$ such that $|Z| = \kappa \subseteq Z$ and $\kappa , \B , A , \dot{A} , a, \gamma \in Z$, we have $Z \cap V_{\gamma} \in S$, hence if $g$ is a $(Z,\B)$-generic filter on $\B$ %(in particular, $Z \cap V_{\gamma}$ is in $S_0$ and hence in $S$), then $Z$ is \lq good', i.e., 
$\dot{A}^g_Z = A \cap M_Z [\bar{g}]$. The rest of the arguments for the tree representation of $A$ as in Theorem~\ref{B-Baire char}~\ref{enum-thm:B-Baire-char3} is exactly the same as in those for the direction \ref{enum-thm:B-Baire-char2}$\Rightarrow$\ref{enum-thm:B-Baire-char3} in Theorem~\ref{B-Baire char}.
\end{description}
This finishes the proof of Theorem~\ref{uB generic absoluteness}.
\end{proof}

We now introduce the notion of universally Baire subsets of $2^{\kappa}$:
\begin{Def}\label{universally Baire}
Let $\Gamma$ be a class of complete Boolean algebras and let $A$ be a subset of $\kappa$. We say $A$ is {\it universally Baire in $2^{\kappa}$ with respect to $\Gamma$} ($\mathsf{uB}^{\Gamma}_{\kappa}$) if $A$ is $\B$-Baire in $2^{\kappa}$ for all $\B$ in $\Gamma$.
\end{Def}

\vspace{1cm}

\begin{Ex}\label{Ex}

${}$

\begin{enumerate}[(1)]
\item Let $\kappa = \omega$ and $\Gamma$ be the class of all complete Boolean algebras. Then $\mathsf{uB}^{\Gamma}_{\kappa}$ sets are exactly those which are universally Baire sets of reals introduced by Feng, Magidor, and Woodin~\cite{MR1233821}. Hence the notion of $\mathsf{uB}^{\Gamma}_{\kappa}$ sets generalizes the classical notion of universally Baire sets of reals.

\item Suppose that the axiom $\FA_{\kappa}(\B)$ {\it fails} for all $\B$ in $\Gamma$. Then {\it every} subset of $2^{\kappa}$ is $\mathsf{uB}^{\Gamma}_{\kappa}$. In fact, by Remark~\ref{meager-FA}, for all $\B$ in $\Gamma$, the whole space $\text{St}(\B)$ is $\kappa$-meager and hence every subset of $2^{\kappa}$ is $\B$-Baire in $2^{\kappa}$. Hence, when we consider whether a given set is $\mathsf{uB}^{\Gamma}_{\kappa}$ or not, we may assume that the axiom $\FA_{\kappa}(\B)$ holds for all $\B$ in $\Gamma$.

%%%% FOR THE ITEM BELOW, CHECK IF B FA_{\KAPPA}(B  RESTRICT b) IS NEEDED FOR ALL b IN B. ALSO, CHECK IF WE NEED PROPER CLASS MANY WOODIN CARDINALS.

\item \label{BPFA} (Caicedo and Velickovic~\cite{MR2231126}) Suppose that there are proper class many Woodin cardinals and assume that $\mathsf{BPFA}$ holds. Let $\kappa = \omega_1$ and $\Gamma$ be the class of all complete Boolean algebras $\B$ such that $\FA_{\omega_1}(\B \upharpoonright b)$ holds for all $b$ in $\B$, and that $\B$ preserves $\mathsf{BPFA}$. Then there is a $\mathsf{uB}^{\Gamma}_{\omega_1}$ well-order on $\wp(\omega_1)$ by the following arguments: By the result of Caicedo and Velickovic~\cite{MR2231126}, under $\mathsf{BPFA}$, there is a well-order on $\wp (\omega_1)$ which is $\Delta_1$ in the structure $(H_{\omega_2}, \in)$ with a parameter of a subset of $\omega_1$. Fix a formula $\psi$ and a parameter $a$ defining the well-order on $\wp (\omega_1)$ over $(H_{\omega_2}, \in)$, and let $\phi (x,a)$ state that $(H_{\omega_2}, \in ) \vDash \psi [x,a]$. Let $\B$ be a complete Boolean algebra in $\Gamma$. We verify the condition \ref{uB-generic-absoluteness2} for $\kappa = \omega_1$, $\phi$, $a$, and $\B$ to conclude that the well-order is $\B$-Baire. Let $G, H , j , M$ be as in the condition \ref{uB-generic-absoluteness2}. Then since $G \in V[H]$ and $M$ is closed under $\omega_1$-sequences in $V[H]$, $(H_{\omega_2} , \in )^{V[G]}$ is a substructure of $(H_{\omega_2}, \in )^M$. Since $\B \in \Gamma$, $V[G]$ is a model. Also, since $V$ is a model of $\mathsf{BPFA}$, by the elementarity of $j\colon V \to M$, $M$ is also a model of $\mathsf{BPFA}$.  Now since $\psi$ is equivalent to a $\Delta_1$ statement in $(H_{\omega_2}, \in)$ under $\mathsf{BPFA}$, the statement $\psi (x, a)$ is absolute between $(H_{\omega_2} , \in )^{V[G]}$ and $(H_{\omega_2}, \in )^M$. Hence for all $x \in \bigl(\wp(\omega_1) \times \wp(\omega_1) \bigr) \cap V[G]$, $V[G] \vDash \phi [x,a] \iff M \vDash \phi [x,a]$. Also, since the critical point of $j$ is $\omega_2^V$ and $a$ is a subset of $\omega_1$ in $V$, $j(a) = a$. Therefore, for all $x\in \bigl(\wp(\omega_1) \times \wp(\omega_1)\bigr) \cap V[G]$, $V[G] \vDash \phi [x,a] \iff M \vDash \phi [x, j(a)]$, as desired.

%This easily follows from the result by Caicedo and Velickovic~\cite{MR2231126} and Theorem~\ref{uB generic absoluteness} in this paper.

\end{enumerate}
\end{Ex}

Similar arguments to Example~\ref{Ex}~\ref{BPFA} give us the following: Suppose there are proper class many Woodin cardinals. Let $\Gamma$ be the class of all complete Boolean algebras $\B$ such that $\FA_{\kappa}(\B \upharpoonright b)$ holds for all $b$ in $\B$. Then if $A \subseteq 2^{\kappa}$ is $\Delta_1$ in the structure $(H_{\kappa^+}, \in )$ in $\mathsf{ZFC}$, then $A$ is $\mathsf{uB}^{\Gamma}_{\kappa}$.

\section{Questions}

We list some open questions regarding universally Baire subsets of $2^{\kappa}$  for an arbitrary $\kappa$ with a special focus in some cases on $\kappa=\omega_1$.

The first question concerns Wadge reducibility for subsets of $2^{\kappa}$: 
Let $A, A'$ be subsets of $2^{\kappa}$. We say $A$ is {\it Wadge reducible to}  $A'$ ($A \le_{\text{W}} A'$) if there is a continuous function $f \colon 2^{\kappa} \to 2^{\kappa}$ such that $A = f^{-1} [A']$.
\begin{Q}
Assume large cardinal assumptions of your preference (e.g., there are proper class many Woodin cardinals). 
Let $\Gamma$ be the class of complete Boolean algebras $\B$ such that for all $b$ in $\B$, 
$\FA_{\kappa} (\B \upharpoonright b)$ holds. 
Let $A, A'$ be $\mathsf{uB}^{\Gamma}_{\kappa}$ subsets of $2^{\kappa}$. 

Can one prove that either $A \le_\text{W} A'$ or $(2^{\kappa} \setminus A') \le_{\text{W}} A$? 
Also, can one prove that the order $\le_{\text{W}}$ on $\mathsf{uB}^{\Gamma}_{\kappa}$ sets is well-founded?
\end{Q}

Note that when $\kappa = \omega$ and $\Gamma$ is the class of all complete Boolean algebras, the answers to both questions are \lq\lq Yes":
this follows from the Wadge Lemma, and the well-foundedness of the Wadge order for universally Baire sets of reals under large cardinals, using the determinacy of Wadge games for universally Baire sets of reals.

The second question is more vague and concerns the possibility to use 
$\mathsf{uB}^{\Gamma}_{\kappa}$ sets to measure the complexity of first-order theories:
\begin{Q}
Let $\Gamma$ be the class of complete Boolean algebras $\B$ such that for all $b$ in $\B$, 
$\FA_{\kappa} (\B \rstrct b)$ holds. Can the theory of $\mathsf{uB}^{\Gamma}_{\kappa}$ 
sets provide tools to measure the model theoretic complexity of mathematical theories?
In particular, is the notion of Borel reducibility in generalized descriptive set theory meaningful to 
compare the complexity of $\mathsf{uB}^{\Gamma}_{\kappa}$ equivalence relations?
\end{Q}

Note that this is the case in descriptive set theory where one uses Borel reducibility to compare complexity of classification problems.

It is known that universally Baire sets of reals are exactly the same as $\infty$-homogeneously  Suslin sets of reals under the existence of large cardinals. 
The third question concerns a characterization of universally Baire subsets of $2^{\kappa}$ in terms of homogeneously Suslin sets:
\begin{Q}
Suppose that there are proper class many Woodin cardinals. 
Let $\Gamma$ be the class of complete Boolean algebras $\B$ such that for all $b$ in $\B$, 
$\FA_{\kappa} (\B \rstrct b)$ holds. 
Can one formulate the definition of homogeneously Suslin subsets of $2^{\kappa}$ which generalizes the case 
$\kappa = \omega$ and then prove that homogeneously Suslin subsets of 
$2^{\kappa}$ are exactly those which are $\mathsf{uB}^{\Gamma}_{\kappa}$ sets?
\end{Q}

It is known that assuming large cardinals, universally Baire sets of reals are exactly those which are generically 
invariant sets of reals (for the precise statement, see e.g.,~\cite[Theorem~3.3.7]{MR2069032}). 
In our context, this is the case when $\kappa = \omega$ and $\Gamma$ is the class of all complete Boolean algebras.

One can ask whether the similar equivalence could be established when  $\kappa = \omega_1$. If $\kappa = \omega_1$, the natural class of forcings one could look at is the class of {\it stationary set preserving} ($\SSP$) forcings, because forcings which are {\it not} $\SSP$ do {\it not} even preserve the $\Sigma_1$-theory of projective subsets of $2^{\omega_1}$. 

In~\cite{MMM}, Viale established a generic absoluteness result for statements about subsets of $2^{\omega_1}$ under the forcing axiom $\MM^{+++}$, which is a natural strengthening of $\MM^{++}$ and of Martin's Maximum.
A forcing axiom of this sort is required, if one aims to obtain the equivalence between $\mathsf{uB}^{\SSP}_{\omega_1}$ sets and $\SSP$-generically invariant sets. 
For the definitions and basics on $\MM^{+++}$, super almost huge cardinals, totally rigid partial orders, and category forcings $\U_{\delta}$, we refer the reader to \cite{MMM}. 
The fourth question is the following:
\begin{Q}
Suppose that the axiom $\MM^{+++}$ holds and that there are proper class many super almost huge cardinals. 
Let $\Gamma$ be the class of $\SSP$-complete Boolean algebras $\B$ which are totally rigid, and force $\MM^{+++}$. 
Is the family of $\mathsf{uB}^{\Gamma}_{\omega_1}$ sets the same as those subsets of $2^{\omega_1}$ which are generically invariant with respect to forcings in $\Gamma$ (i.e., those $A\subseteq 2^{\omega_1}$ defined by a formula $\phi$ which is absolute among $V$, $V^{\B}$ for $\B \in \Gamma$, and generic ultrapowers $M$ obtained by the category forcings $\mathsf{U}_{\delta}\rstrct\B$ where $\delta>|\B|$ is a super almost huge cardinal)?
\end{Q}

%To be filled in.

\bibliographystyle{plain}
\bibliography{myreference}

\end{document}